\documentclass{article}

\usepackage[english]{babel}

\usepackage[
backend=biber,
style=alphabetic,
]{biblatex}

\usepackage[a4paper,top=2cm,bottom=2cm,left=3cm,right=3cm,marginparwidth=1.75cm]{geometry}

\usepackage{enumitem}
\usepackage{comment}
\usepackage{amsmath}
\usepackage{amssymb}
\usepackage{amsfonts}
\usepackage{amsthm}
\usepackage{graphicx}
\usepackage{csquotes}
\usepackage{mathtools}
\usepackage[colorlinks=true, allcolors=blue]{hyperref}
\usepackage[capitalize]{cleveref}
\usepackage{xcolor}

\theoremstyle{plain}
\newtheorem{theorem}{Theorem}
\newtheorem*{theorem*}{Theorem}
\newtheorem{lemma}[theorem]{Lemma}
\newtheorem*{lemma*}{Lemma}
\newtheorem{proposition}[theorem]{Proposition}

\newtheorem{corollary}[theorem]{Corollary}

\theoremstyle{remark}
\newtheorem{remark}[theorem]{Remark}

\DeclarePairedDelimiter\ceil{\lceil}{\rceil}

\newcommand{\length}[1]{\ell(#1)}
\newcommand{\abs}[1]{|#1|}

\newcommand{\eps}{\varepsilon}

\newcommand{\union}{\cup}

\newcommand{\reals}{\mathbb{R}}

\newcommand{\vol}{\mathrm{Vol}}

\newcommand{\mheat}[4]{p_{#1}^{#4}(#2,#3)}

\newcommand{\thin}[1]{{#1}_{\mathrm{thin}}}
\newcommand{\thick}[1]{{#1}_{\mathrm{thick}}}
\newcommand{\truncated}[1]{\tilde{#1}}

\newcommand{\re}[1]{\mathrm{Re} \left(#1\right)}
\newcommand{\inj}[1]{\mathrm{Inj}\left(#1\right)}

\newcommand{\hyperbolic}{{\mathbb{H}}}
\newcommand{\dist}[2]{\mathrm{dist}\left(#1, #2 \right)}
\newcommand{\mdist}[3]{\mathrm{dist}_{#3}\left(#1, #2 \right)}
\newcommand{\invlengths}{L}
\newcommand{\tracediff}[2]{D_{#2}\left({#1}\right)}
\newcommand{\gsmall}{g_{\mathrm{small}}^k}

\newcommand{\one}{\mathbf{1}}
\renewcommand{\S}{{S}}   

\usepackage{soul}

\title{Determinants of Laplacians on converging hyperbolic surfaces}
\author{Renan Gross\footnote{Department of Pure Mathematics and Mathematical Statistics, University of Cambridge. rg751@cam.ac.uk}, Guy Lachman\footnote{Department of mathematical sciences, Tel Aviv University. guy.lachman@gmail.com}, Asaf Nachmias\footnote{Department of mathematical sciences, Tel Aviv University. asafnach@tauex.tau.ac.il}}

\begin{document}
\maketitle

\begin{abstract}
Let $\S_k$ be a sequence of compact hyperbolic surfaces of increasing volume which locally converges to a random rooted surface. We show that if the normalized sum of the reciprocal lengths of very short simple closed geodesics converges to 0, then the normalized logarithm of the determinant of the Laplacian of $\S_k$ converges to a constant depending only the law of the limiting surface. 
\end{abstract}

\section{Introduction}
Let $\S$ be a compact hyperbolic surface, that is, a compact Riemannian surface of constant curvature $-1$. Our main object of study is the \emph{log determinant of the Laplacian}. This is a very classical object so our description here is brief; see e.g. \cite[][Section 1]{awonusika_determinants_of_laplacians} and the references therein for a thorough introduction. The Laplacian $\Delta_{\S}$ has a discrete spectrum $0=\lambda_0<\lambda_1\leq \lambda_2\leq \ldots$ increasing to infinity. Analogously to finite matrices, one would intuitively wish to define the log of the determinant of the Laplacian as the sum of the logarithms of the eigenvalues $\lambda_j$. This is impossible to do directly, as the eigenvalues tend to infinity; rather, we instead use the Laplacian's zeta function, defined as the formal sum
\begin{equation*}
    \zeta_{\S}(s) := \sum_{j=1}^\infty \frac{1}{\lambda_j^s} \, .
\end{equation*}
Using Weyl's law, $\zeta_{\S}(s)$ is well-defined when $\re s$ is large enough, and using the Minakshisundaram-Pleijel asymptotic expansion, it can be analytically extended to a meromorphic function on $\mathbb{C}$ which is analytic at $0$ (see e.g. \cite[][eq 45, chapter VI]{chavel_book}). We then define the \emph{log determinant} of the Laplacian as
\begin{equation*}
    \log \det (\Delta_{\S}) := -\zeta_{\S}'(0) \,
\end{equation*} 
(when treating $\zeta_{\S}(s)$ as a formal series, the expression $-\zeta_{\S}'(0)$ is indeed the sum of the logarithms of the non-zero eigenvalues).

Let $\mathcal{H}$ be the space of all rooted hyperbolic surfaces equipped with the Gromov-Hausdorff metric, and let $(\mathbf{\S}_\infty, \mathbf{x}_0)$ be a random rooted surface in $\mathcal{H}$ with law $\mu$. A compact hyperbolic surface $\S$ defines a measure $\mu_\S$ on $\mathcal{H}$ by pushing forward the uniform measure $\vol_\S(\cdot) / \vol_\S(\S)$ by the map $x \mapsto (\S, x)$. A sequence $\S_k$ of compact hyperbolic manifolds is said to \emph{converge locally} to $(\mathbf{\S}_\infty, \mathbf{x}_0)$ if $\mu_{\S_k}$ weakly converges to $\mu$ as $k \to \infty$ (this is also known as \emph{Benjamini-Schramm} convergence).

\begin{figure}[ht]
    \centering
    \includegraphics[scale = 0.25]{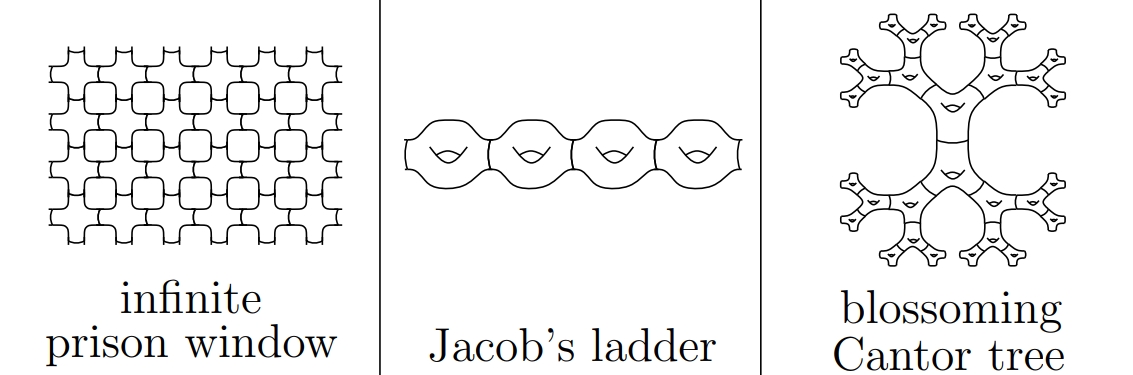}
    \caption{ In addition to the hyperbolic plane, a sequence of hyperbolic surfaces can converge to various infinite structures, such as a thick grid, Jacob's ladder, or a Cantor tree. Image taken from \cite{angel_hutchcroft_nachmias_ray_unimodular_random_maps}.}
    \label{fig:types_of_limiting_shapes}
\end{figure}

The simplest type of convergence of hyperbolic surfaces is to the hyperbolic plane $\hyperbolic$. In this case, the above condition means that for any fixed $R > 0$, a ball of radius $R$ around a uniformly random point in $\S_k$ is isometric to a ball in $\hyperbolic$ of radius $R$ with probability tending to $1$ as $k\to\infty$ (in other words, the injectivity radius seen from a uniformly drawn point tends to $\infty$ with high probability). Convergence to $\hyperbolic$ is typical in the sense that various models of random hyperbolic surfaces have this property, see \cite{naud_determinants_of_laplacians}. However, convergence to non-simply connected infinite rooted surfaces is also possible (e.g. hyperbolic structures on the topologies depicted in Figure \ref{fig:types_of_limiting_shapes}).

In this paper, inspired by the work of Naud \cite{naud_determinants_of_laplacians}, we provide a sharp condition which implies that the log determinant of $\S_k$, when properly normalized, converges to a constant depending only on the distribution of the limit surface $(\mathbf{\S}_\infty, \mathbf{x}_0)$. For $\eta >0$, let $\{\ell_i\}_{i=1}^{s}$ be the lengths of all simple closed geodesics in $\S_k$ of length at most $\eta$, and define
\begin{equation*}
    \invlengths_\eta (\S) = \frac{1}{\vol(\S)} \sum_{i=1}^{s} \frac{1}{\ell_i} \, .
\end{equation*}
We say that a sequence of surfaces $\S_k$ has \emph{uniformly integrable geodesics} if for every $\eps > 0$, there exists an $\eta > 0$ such that $\invlengths_\eta(\S_k) < \eps$ for all $k$ large enough. 

\begin{theorem}\label{thm:logdet_main_theorem}
    There exists a universal constant $E_\hyperbolic > 0$ so that the following holds. Let $\S_k$ be a sequence of compact hyperbolic surfaces of volume tending to $\infty$ that locally converges to a random rooted surface $(\mathbf{\S}_\infty, \mathbf{x}_0)$ with distribution $\mu$. Define $E_{\mu}\in [-\infty,\infty)$ by
    \begin{equation} \label{eq:main_constant_to_converge_to}
        E_{\mu} = E_\hyperbolic - \int_0^\infty \frac{\int [\mheat{t}{x_0}{x_0}{\S_\infty} - \mheat{t}{0}{0}{\hyperbolic}] d\mu}{t} ~dt \, ,
    \end{equation}
    where $\mheat{t}{x}{y}{\S}$ is the heat-kernel of the Laplacian on a Riemannian surface $\S$.
    \begin{enumerate}
        \item \label{enu:main_theorem_uniform_integrable} If $\S_k$ has uniformly integrable geodesics, then
        \begin{equation*}
            \lim_{k\to\infty} \frac{\log \det \Delta_{\S_k}}{\vol(\S_k)} = E_{\mu} \, .
        \end{equation*}

        \item \label{enu:main_theorem_non_uniform_integrable} If $\S_k$ does not have uniformly integrable geodesics, then
        \begin{equation*}
            \limsup_{k\to\infty} \frac{\log \det \Delta_{\S_k}}{\vol(\S_k)} \leq E_\mu \, ,
        \end{equation*}
        with strict inequality if the limit exists and is finite. 
    \end{enumerate}    
\end{theorem}
Note that if $\S_k$ converges to the hyperbolic plane, the numerator in the integrand on the right-hand side of \eqref{eq:main_constant_to_converge_to} is identically $0$, and so the log determinant simply converges to $E_\hyperbolic$. The converse is also true.
\begin{corollary}\label{cor:convergence_to_plane}
    Let $\S_k$ be a sequence of compact hyperbolic surfaces of volume tending to $\infty$ that locally converges to a random rooted surface $(\mathbf{\S}_\infty, \mathbf{x}_0)$ with distribution $\mu$. Assume that $\S_k$ has uniformly integrable geodesics and that 
    \begin{equation*}
        \lim_{k\to\infty} \frac{\log \det \Delta_{\S_k}}{\vol(\S_k)} = E_{\hyperbolic} \, .
    \end{equation*}
    Then $\S_k$ converges to the hyperbolic plane, i.e. $\mathbf{\S_\infty} = \hyperbolic$ $\mu$-a.s.
\end{corollary}
The proof of \cref{thm:logdet_main_theorem} relies on bounding both the divergence of the heat kernel $\mheat{t}{x}{x}{\S_k}$ at small times $t \to 0$ as well as the difference with its limiting value of $1/\vol(\S_k)$ at large times $t\to \infty$. Indeed, for a compact hyperbolic surface $\S$, if we denote
\begin{equation*}
    \tracediff{t}{\S} = \int_{\S} \big [\mheat{t}{x}{x}{\S} - \mheat{t}{0}{0}{\hyperbolic}\big ] ~dx \, ,
\end{equation*}
then a  calculation involving $\zeta(s)$ and the heat trace formula (see \cite[][section 2]{naud_determinants_of_laplacians}) gives that
\begin{equation}\label{eq:the_starting_point}
    \log \det \Delta_{\S} = \vol(\S)E_\hyperbolic + \gamma_0 - \int_0^1\frac{\tracediff{t}{\S}}{t}~dt - \int_1^\infty \frac{\tracediff{t}{\S}-1}{t} ~dt \, ,
\end{equation}
where $\gamma_0$ is the Euler-Mascheroni constant. This is the starting point of the proof of \cref{thm:logdet_main_theorem}. We will treat each of the two terms $\frac{1}{\vol(\S_k)}\int_0^1\frac{\tracediff{t}{\S_k}}{t}~dt$ and $\frac{1}{\vol(\S_k)} \int_1^\infty \frac{\tracediff{t}{\S_k}-1}{t}~dt$ separately. For large times, we provide a bound depending only on the sum of reciprocal lengths of geodesics:
\begin{lemma} \label{lem:ds_bounded_large_times}
    There exists a constant $C > 0$ such that for every compact hyperbolic surface $\S$ and every $t>0$,
    \begin{equation*}
        \frac{\abs{\tracediff{t}{\S}-1}}{t~\vol(\S)} \leq C \left( \frac{e^{-t/4}}{t^2} + \frac{\sqrt{1 + \invlengths_{2\sinh^{-1}(1)}(\S)}}{t^{3/2}} \right) \, .
    \end{equation*}    
\end{lemma}
For small times, we recall the definition of a uniformly integrable family of functions (from which the term \textit{uniformly integrable geodesics} is derived): Let $\mu$ be a finite measure on some space $X$. A family of integrable functions $f_k : X \to \reals$ is called \emph{uniformly integrable} if they are uniformly bounded in $\mathrm{L}^1(X)$, and if for all $\eps > 0$, there exists a $\delta > 0$ such that if $A \subseteq X$ has measure $\mu(A) \leq \delta$, then $\sup_k \int_A \abs{f_k(x)}~d\mu(x) \leq \eps$. 
\begin{lemma} \label{lem:ds_uniformly_integrable_small_times}
    The family of functions $\Big\{ \frac{\tracediff{t}{\S_k}}{t~\vol(\S_k)}\Big\}_{k=1}^{\infty}$ is uniformly integrable on $(0,1)$ if and only if the sequence $\S_k$ has uniformly integrable geodesics.
\end{lemma}

In the proof we will use two external ingredients. The first is a recent result of the authors bounding the convergence rate of $\mheat{t}{x}{x}{\S}$ to the uniform measure. 

    \begin{theorem}[\cite{gross_lachman_nachmias_sharp_lower_bounds}, Theorem 2 and Corollary 9]\label{thm:convergence_to_uniform_large_t}There exists a constant $C > 0$ such that for any compact hyperbolic surface $\S$,
    \begin{equation*}
         \frac{1}{\vol(\S)} \int_\S \Big| \mheat{t}{x}{x}{\S} -\frac{1}{\vol(\S)} \Big | ~dx \leq C~ \sqrt{\frac{1 + \invlengths_{2\sinh^{-1}(1)}(\S)}{t}} \qquad \forall t \geq 1 \, . 
    \end{equation*}
    \end{theorem}

The second is a classical estimate of Li and Yau \cite{li_yau_main}. 

\begin{theorem}[Bounding the heat kernel by volume of balls; Corollary 3.1 in \cite{li_yau_main} with $\alpha=3/2$ and $\varepsilon=1/2$] \label{thm:li_yau}
    Let $\S$ be a complete Riemannian manifold without boundary and with Ricci curvature bounded from below by $-K$. Then the heat-kernel on $\S$ satisfies
    \begin{equation*} 
        \mheat{t}{x}{y}{\S} \leq C_1\vol(B(x,\sqrt{t}))^{-1/2}\vol(B(y,\sqrt{t}))^{-1/2}\exp\left(C_2 Kt-\frac{d_{\S}(x,y)^2}{4.5t} \right) \, .
    \end{equation*}
    \end{theorem}

\subsection{Previous work and remarks} \label{subsec:previous_work}
Our result strengthens a recent paper by Naud \cite{naud_determinants_of_laplacians} by eliminating the need for a spectral gap and loosening the condition on short closed geodesics; further, it allows $\S_k$ to locally converge to any hyperbolic surface, not just to the hyperbolic plane (Naud's conditions implicitly imply convergence to $\hyperbolic$). A similar convergence result has been proved for cyclic covers of general Riemannian manifolds \cite{dang_lin_naud_asymptotics_zeta_determinants}.

The discrete counterpart to the determinant of the Laplacian is the quantity known as \emph{tree entropy} defined by Lyons \cite{lyons_asymptotic_enumeration_of_spanning_trees} (the relation to trees is not surprising, indeed Kirchhoff's theorem relates the product of the non-zero eigenvalues of a graph's Laplacian with the number of its spanning trees). Our main theorem (\cref{thm:logdet_main_theorem}) is the continuous analogue of \cite[Theorem 3.2]{lyons_asymptotic_enumeration_of_spanning_trees} and our approach in this paper is partially inspired by that theorem. 

In \cite{abert_et_al_growth_of_l2_invariants}, it is shown that other properties of converging manifolds carry over to the limit (in both two and higher dimensions): Betti numbers (\cite[Theorem 7.13]{abert_et_al_growth_of_l2_invariants}) and analytic torsion (\cite[Theorem 8.4]{abert_et_al_growth_of_l2_invariants}, under bounded injectivity conditions and acyclicity conditions). These results do not straightforwardly imply \cref{thm:logdet_main_theorem}.

\begin{remark}
    The constant $E_\mu$ can be morally thought of as the derivative of the normalized zeta function of the distribution $\mu$, as follows. For a compact manifold $\S$, the zeta function is given by \cite[][eq. (1.3)]{sarnak_determinants_of_laplacians_heights}
    \begin{align*}
        \zeta_{\S}(s) &= \frac{1}{\Gamma(s)} \int_{0}^{\infty} t^{s-1}\left(\mathrm{Tr}\left(e^{-t\Delta_{\S}}\right) -1 \right)dt \\
        &= \frac{1}{\Gamma(s)} \int_{0}^{\infty} t^{s-1}\left(\int_{\S} \mheat{t}{x}{x}{\S} ~dx -1 \right)dt \,.
    \end{align*}
    When normalizing by the volume of the manifold, denoting by $\mu_{\S}$ the uniform distribution on $\S$, we have 
    \begin{equation*}
        \frac{\zeta(s)}{\vol(\S)} = \frac{1}{\Gamma(s)} \int_{0}^{\infty} t^{s-1}\left(\int \mheat{t}{x}{x}{\S}d\mu_{\S} -\frac{1}{\vol(\S)} \right)dt \, .
    \end{equation*}
    When $(\mathbf{\S}_\infty, \mathbf{x}_0)$ is a random, rooted, infinite-volume surface with distribution $\mu$, we can analogously define
    \begin{equation*}
        \zeta_\mu(s) := \frac{1}{\Gamma(s)} \int_{0}^{\infty} t^{s-1}\left(\int \mheat{t}{\mathbf{x}_0}{\mathbf{x}_0}{\mathbf{\S}_\infty}d\mu \right)dt \, .
    \end{equation*}
    Under certain conditions (e.g. when $\mathbf{\S}_\infty$ has a spectral gap and bounded injectivity radius a.s.), the constant $E_\mu$ can then be shown to equal 
    \begin{equation*}
        E_\mu = -\zeta_{\mu(s)}'(0) \,.
    \end{equation*}
\end{remark}

\begin{remark} 
    Let $\S$ be a compact hyperbolic surface. In \cite[][Corollary 9]{gross_lachman_nachmias_sharp_lower_bounds}, it was shown that if $C$ is a hyperbolic collar whose waist is a geodesic of length $\ell$, then 
    \begin{equation*}
        \int_C \frac{1}{\inj{x,\S}^2 \land 1} ~dx
    \end{equation*}
    is proportional to $1/\ell$, where $\inj{x,\S}$ is the injectivity radius of $x$. Thus, having uniformly bounded geodesics is equivalent to the family of functions $\Big\{\frac{1}{\inj{x,\S_k}^2 \land 1} \Big\}$ being uniformly integrable.
\end{remark}

\section{Injectivity radius and heat kernel bounds}
Throughout this paper we use the term  ``universal constant'' to mean a real number that does not depend on $\S$, $k$, or any other parameter. For two non-negative functions $f(x)$ and $g(x)$, we write $f \lesssim g$ to mean that there exists a universal constant $C > 0$, such that $f(x) \leq C g(x)$ for all $x$. We write $f(x) \asymp g(x)$ if both $f \lesssim g$ and $g \lesssim f$. 

For a point $x \in \S$, we write $B(x,r)$ for the closed ball in $\S$ of radius $r$ around $x$. The injectivity radius of $x \in \S$, denoted $\inj{x}$, is the supremum radius $r$ such that $B(x,r)$ isometric to the ball of radius $r$ in the hyperbolic plane $\hyperbolic$.

An important geometric characterization of compact hyperbolic surfaces is the so-called collar lemma, which describes the parts of the surface with small injectivity radius. The following is a combination of Theorem 4.1.1 and Theorem 4.1.6 from \cite{buser_book}.
\begin{lemma}[Collar lemma] \label{lem:collar_lemma}
    Let $\gamma_1, \ldots, \gamma_s$ be the set of all simple closed geodesics of length $\leq 2\sinh^{-1}(1)$ on a hyperbolic surface $\S$. Let $W(\gamma_i) =\sinh^{-1}\left(\frac{1}{\sinh\left(\frac{1}{2} \length{\gamma_i}\right)}\right)$ and $C(\gamma_i) = \{x\in \S \mid \mathrm{dist}(x,\gamma_i) \leq W(\gamma_i)\}$. Then 
    \begin{itemize}
        \item The sets $C(\gamma_i)_{i=1}^s$ are pairwise disjoint.
        \item $\inj{x,\S} \geq \sinh^{-1}(1)$ for all $x \notin \union_i C(\gamma_i)$.
        \item Each $C(\gamma_i)$ is isometric to the cylinder $\left[-W(\gamma_i), W(\gamma_i)\right] \times S^1$ with the Riemannian metric $ds^2 = d\rho^2 + \length{\gamma_i}^2\cosh(\rho)^2 d\theta^2$.
        \item If $x \in C(\gamma_i)$ is such that $\inj{x} \leq \sinh^{-1}(1)$ and $d = \mdist{x}{\partial C(\gamma_i)}{\S}$, then
        \begin{equation*}
            \sinh(\inj{x}) = \cosh\left(\frac{1}{2}\length{\gamma_i}\right)\cosh(d) - \sinh(d).
        \end{equation*}
    \end{itemize}
\end{lemma}
Although the lemma gives an exact expression for the injectivity radius of points inside collars as a function of their distance from the boundary, we will mostly make use of the following simple approximation (see e.g. \cite[][Proposition 8]{gross_lachman_nachmias_sharp_lower_bounds} for a proof).

\begin{proposition}[Injectivity radius estimate]\label{prop:injradius_estimate}
    Let $T$ be an infinite hyperbolic collar whose shortest closed geodesic $\gamma$ has length $\ell$. Denote by $r(\rho,\ell)$ the injectivity radius at a point of distance $\rho$ from the waist. Then
    \begin{equation}\label{eq:alternate_injectivity_in_collar}
        r(\rho, \ell) = \frac{1}{2}\cosh^{-1}\left(1 + \left(\cosh(\ell) -1 \right)\cosh(\rho)^2 \right) \, .
    \end{equation}
    In particular, 
    \begin{equation} \label{eq:injradius_estimate}
        r(\rho, \ell) \asymp \ell \cosh(\rho) \, .
    \end{equation}
    for all $\ell \leq 2\sinh^{-1}(1)$ and all $\rho \leq W(\gamma)$ where $W(\gamma)$ is defined in \cref{lem:collar_lemma}.
\end{proposition}

The following lemma gives a bound on the difference between the heat kernels of a compact surface and the hyperbolic plane. 

\begin{lemma} \label{lem:counting_argument_by_injradius}
    For every compact hyperbolic surface $\S$, every $x \in \S$ and every $t>0$, we have both
    \begin{equation} \label{eq:ft_bound}
        \mheat{t}{x}{x}{\S} - \mheat{t}{0}{0}{\hyperbolic} \lesssim \left(1+\frac{1}{\inj{x}}\right) \left(\frac{1}{\sqrt{t}} + \sqrt{t} \cdot e^{64t} \right) \, 
    \end{equation}
    and 
    \begin{equation} \label{eq:sqrt_time_bound}
        \mheat{t}{x}{x}{\S} - \mheat{t}{0}{0}{\hyperbolic} \lesssim \left(1 + \frac{1}{\inj{x}^3} \right) \sqrt{t} \cdot e^{64t}\,.
    \end{equation}
\end{lemma}
The bound \eqref{eq:ft_bound} is optimal in terms of its dependence on small values of the injectivity radius; however, the bound grows to infinity as $t \to 0$. On the other hand, the bound \eqref{eq:sqrt_time_bound} vanishes as $t\to 0$, at the expense of a worse dependence on the injectivity radius. This can be seen as an explicit error bound in the the Minakshisundaram-Pleijel asymptotic expansion. The bound \eqref{eq:ft_bound} is already known (see e.g. Lemma 7.11 of \cite{abert_et_al_growth_of_l2_invariants}); since the proof of both bounds share a common initial calculation, we give a self contained elementary proof of both.

\begin{remark}
    The powers of $t$ and $\inj{x}$ cannot be separately improved in \eqref{eq:sqrt_time_bound}. For example, when $t = \ell^2$ and for points $x$ near the waist of a collar corresponding to a closed geodesic of length $\ell$, it can be seen that $\mheat{t}{x}{x}{\S} - \mheat{t}{0}{0}{\hyperbolic} \gtrsim \frac{1}{\ell^2}$.
\end{remark}

\begin{proof}
    Recall that the heat kernel $\mheat{t}{x}{y}{\S}$ of any hyperbolic surface $\S$ is given by $\mheat{t}{x}{y}{\S} = \sum_{T \in \Gamma} \mheat{t}{\Tilde{x}}{T\Tilde{y}}{\hyperbolic}$, where $\Gamma$ is a group of isometries acting on $\hyperbolic$ so that $\S = \Gamma \backslash \hyperbolic$, and $\Tilde{x}$ and $\Tilde{y}$ are arbitrary inverse images of $x$ and $y$ of the projection to $\S$. To estimate $\mheat{t}{x}{x}{\S}$, we partition the group translations by distance of $T \tilde x$ from $\tilde x$.  Denoting
    \begin{equation*}
        \Gamma(m) := \{T \in \Gamma \mid m < d_\hyperbolic(\tilde x, T \tilde x) \leq m+1 \} \, ,
    \end{equation*} 
    we have
    \begin{align}
       \mheat{t}{x}{x}{\S} - \mheat{t}{0}{0}{\hyperbolic} &=  \sum_{m=0}^{\infty} \sum_{T \in \Gamma(m)} \mheat{t}{\tilde x}{T\tilde x}{\hyperbolic} \nonumber \\ 
       &=\sum_{T \in \Gamma(0)} \mheat{t}{\tilde x}{T \tilde x}{\hyperbolic}\ +  \sum_{m=1}^{\infty} \sum_{T \in \Gamma(m)} \mheat{t}{\tilde x}{T \tilde x}{\hyperbolic} \, . \label{eq:partition_of_the_count}  
    \end{align}
    Let $r = \inj{x}$. The number of group elements $T \in \Gamma(m)$ can be estimated using a simple volume argument (see \cite[][Lemma 7.5.3]{buser_book}): when $D$ is a disc of radius $r/2$ around $\tilde{x}$, the translates $T(D)$ are disjoint for different $T \in \Gamma(m)$. Since the volume of a hyperbolic disc of radius $R$ is $2\pi (\cosh(R)-1)$, this argument yields
    \begin{equation*}
        \#\Gamma(m) \leq \frac{\cosh(m+1+r/2)-1}{\cosh(r/2)-1} \leq \frac{e^{r/2+1}}{\cosh(r/2)-1} e^m \lesssim \left(1+\frac{1}{r^2} \right) e^m \, ,
    \end{equation*}
    where in the inequalities we relied on the asymptotics of the function $\cosh(r)$ for $r \to 0$ and $r \to \infty$.    
    
    A better bound for $x$ with small injectivity radius may be obtained by counting the number of group elements which send $\tilde{x}$ to only a small distance away, as follows.

    \begin{figure}
        \centering
        \includegraphics[width=0.75 \textwidth]{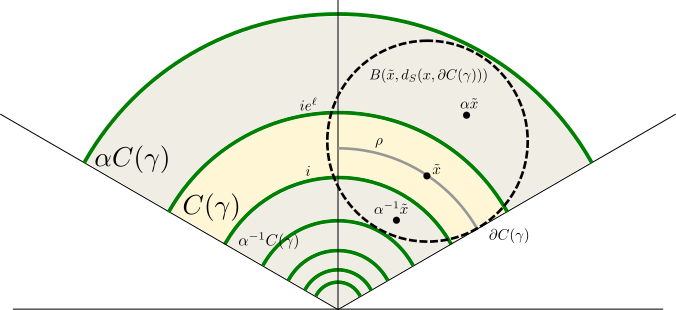}
        \caption{When $\tilde x$ is far from the boundary of the the collar, the group elements which move it the least are the collar translations.} 
        \label{fig:tiling_cone_by_collars}
    \end{figure}

    Assume first that $r \leq \frac{1}{8}\sinh^{-1}(1)$. By \cref{lem:collar_lemma}, $x$ is inside a collar $C(\gamma)$ of width $W(\gamma)$, where $\gamma$ is a simple closed geodesic. In fact, $x$ is at distance strictly greater than $1$ from the boundary $\partial C$ of the collar. This means that all elements $T \in \Gamma(0)$ are of the form $\alpha^k$, where $\alpha$ is the element which generates the collar and $k$ is a nonzero integer. To see this, consider Figure \ref{fig:tiling_cone_by_collars}. The area in yellow is a lift of the collar $C(\gamma)$ to the upper-half plane, given by all points at  hyperbolic distance at most $W(\gamma)$ from the interval $[i, ie^{\length{\gamma}}]$. The translates of this region by powers of $\alpha$ tile a cone whose tip is at the origin. Any element in $\Gamma$ which is not a power of $\alpha$ must take $\tilde x$ outside this cone. Using elementary hyperbolic geometry, it can be shown that the ball of radius $\mdist{x}{\partial C(\gamma)}{\S}$ around $\tilde x$ is tangent to the cone's boundary, which means that $\mdist{\tilde x}{y}{\hyperbolic} > 1$ for all $y$ outside the cone, i.e. we necessarily have $\mdist{\tilde x}{T \tilde x}{\hyperbolic} > 1$ for all $T \neq \alpha^k$. 
    
    Denote the hyperbolic distance from $\tilde{x}$ to the imaginary axis by $\rho$. The distance $\mdist{\tilde x}{\alpha^k \tilde x}{\hyperbolic}$ is then equal to twice the injectivity radius of a point at distance $\rho$ from the central geodesic of an infinite collar with waist $k\ell(\gamma)$. By equation \eqref{eq:alternate_injectivity_in_collar}, we have
    \begin{equation} \label{eq:distance_of_translates}
        \mdist{\tilde x}{\alpha^k \tilde x}{\hyperbolic} = \cosh^{-1}\left(1+\left(\cosh (k \length{\gamma}) - 1  \right) \cosh (\rho)^2 \right) \, .
    \end{equation}
    Every $k$ such that $\alpha^k \in \Gamma(0)$ thus satisfies
    \begin{equation*}
        \cosh(1) \geq 1 + (\cosh(k\ell(\gamma)) -1) \cosh(\rho)^2 \, .
    \end{equation*}
    Rearranging and using the fact that $\cosh(z) \geq 1 + \frac{z^2}{2}$, this gives 
    \begin{equation*}
        \abs{k} \leq \sqrt{\frac{2(\cosh(1)-1)}{\cosh(\rho)^2 \ell^2}} \asymp \frac{1}{\inj{x}} \, ,
    \end{equation*}
    where the last relation is by \eqref{eq:injradius_estimate}. Thus, there exists a universal constant $C > 0$ such that 
    \begin{equation} \label{eq:counting_closeby_points}
        \#\Gamma(0) \leq  \frac{C}{\inj{x}} \, .
    \end{equation}
    Let $D$ be the disc of radius $1/2$ around $\tilde{x}$, and consider the translates $T(D)$ for $T \in \Gamma(m)$. These discs are contained in $B(\tilde{x}, m+2)$, but need not be disjoint. However, no point $y \in B(\tilde{x}, m+2)$ can be contained in too many such discs. To see this, let $T_1,\ldots, T_s \in \Gamma(m)$ be isometries such that $y \in B(T_i \tilde{x}, 1/2)$ for all $i = 1, \ldots, s$. Then the disc $B(T_1 \tilde{x}, 1)$ contains all the points $T_i\tilde{x}$. The image of this disc under $T_1^{-1}$ is a disc of radius $1$ around $\tilde{x}$, containing $s$ distinct points $T_1^{-1}T_i\tilde{x}$. Each isometry $T_1^{-1}T_i$ ($i\neq 1$) is therefore a member of $\Gamma(0)$, and by \eqref{eq:counting_closeby_points}, we must have $s-1 \leq C/\inj{x}$. The sum of volumes of $T(D)$ over $T\in \Gamma(m)$ can thus be no more than $1+C / \inj{x}$ times the volume of $B(\tilde{x}, m+2)$, yielding
    \begin{equation*}
        \#\Gamma(m) \lesssim \frac{1}{\inj{x}}e^m \,.
    \end{equation*}
    Finally, relaxing the demand that $r \leq \frac{1}{8}\sin^{-1}(1)$, we have the general result that 
    \begin{equation}\label{eq:counting_group_actions}
        \#\Gamma(m) \lesssim \left(1 + \frac{1}{\inj{x}}\right) e^m \,.
    \end{equation}

    We now return to bounding \eqref{eq:partition_of_the_count}. A well-known bound on $p_t^{\hyperbolic}$ (see \cite[][Lemma 7.4.26]{buser_book}) states that
    \begin{equation} \label{eq:simple_hyperbolic_heat_kernel_bound}
        \mheat{t}{x}{y}{\hyperbolic} \lesssim \frac{1}{t}\exp\left(-{\frac{{d_\hyperbolic(x,y)}^2}{8t}}\right) \, .
    \end{equation}
    The second term on the right-hand side of \eqref{eq:partition_of_the_count} can then be bounded by
    \begin{equation*}
        \sum_{m=1}^{\infty} \sum_{T \in \Gamma(m)} \mheat{t}{\tilde x}{T \tilde x}{\hyperbolic} \lesssim \left( 1+ \frac{1}{r}\right) \frac{1}{t} \sum_{m=1}^\infty \exp\left( m - m^2 / 8t \right) \, .
    \end{equation*}
    Since $\frac{1}{x}\exp(-1/x) \leq 1$ for all $x$, we have that $\frac{1}{t}\exp(- m^2/8t) \lesssim \exp(- m^2/16t)$ for all $m \geq 1$, giving 
    \begin{equation*}
        \sum_{m=1}^{\infty} \sum_{T \in \Gamma(m)} \mheat{t}{\tilde x}{T \tilde x}{\hyperbolic} \lesssim \left( 1+ \frac{1}{r}\right) \sum_{m=1}^\infty \exp\left( m - m^2 / 16t \right) \, .
    \end{equation*}
    Whenever $m \geq 32t$, we have $m \leq m^2/32t$, and so the sum in the above equation can be bounded from above by
    \begin{align}
        \sum_{m=1}^\infty \exp\left( m - m^2 / 16t \right) &\leq \sum_{m=1}^{\ceil{32t}} \exp\left( m - m^2 / 16t \right) + \sum_{m=\ceil{32t}}^\infty \exp\left( - m^2 / 32t \right) \nonumber \\
        &\leq \ceil{32t} e^{\ceil{32t}} \exp(-1/16t) + \sum_{m=\ceil{32t}}^\infty \exp\left( - m^2 / 32t \right) \, . \label{eq:counting_mid2_calculation}
    \end{align}
    Since the function $m \mapsto \exp\left(-m^2/32t\right)$ is decreasing, we can bound the second sum on the right-hand side via a Gaussian integral bound:
    \begin{equation*}
        \sum_{m=\ceil{32t}}^\infty \exp\left( - m^2 / 32t \right) \leq \sum_{m=1}^\infty \exp\left( - m^2 / 32t \right) \leq \int_0^\infty \exp\left(-x^2/32t\right) ~dx \lesssim \sqrt{t} \, .
    \end{equation*}
    Using the simple bound $xe^x \leq \sqrt{x}e^{2x}$ to bound the first term on the right-hand side in \eqref{eq:counting_mid2_calculation}, we get
    \begin{equation} \label{eq:large_translations}
        \sum_{m=1}^{\infty} \sum_{T \in \Gamma(m)} \mheat{t}{\tilde x}{T \tilde x}{\hyperbolic} \lesssim \left(1 + \frac{1}{r} \right)  \sqrt{t} \cdot e^{64t} \, ,
    \end{equation}
    which is bounded above by the right-hand side of both \eqref{eq:ft_bound} and \eqref{eq:sqrt_time_bound}.
    
    It remains to bound the first term on the right-hand side of \eqref{eq:partition_of_the_count}. If $r > \frac{1}{8}\sinh^{-1}(1)$, we use the fact that $d_{\hyperbolic}(\tilde{x},T\tilde{x}) \geq r$ for every $T \neq \mathrm{Id}$; by \eqref{eq:counting_group_actions} and \eqref{eq:simple_hyperbolic_heat_kernel_bound} we have
    \begin{equation}
        \sum_{T \in \Gamma(0)} \mheat{t}{\tilde x}{T \tilde x}{\hyperbolic} 
        \lesssim \left( 1 + \frac{1}{r}\right) \exp \left( -r^2 / 8t\right) \lesssim \left( 1 + \frac{1}{r}\right) \sqrt{t} \, , \label{eq:small_translations_large_injradius}
    \end{equation}
    where the last inequality is since $r$ is bounded from below. This again is  bounded above by the right-hand side of both \eqref{eq:ft_bound} and  \eqref{eq:sqrt_time_bound}.

    Suppose now that $r \leq \frac{1}{8}\sinh^{-1}(1)$. As we saw above, in this case $x$ is in a collar $C(\gamma)$, and all elements $T \in \Gamma(0)$ are of the form $\alpha^k$, where $\alpha$ is the element which generates the collar and $k$ is a non-zero integer. We therefore have
    \begin{equation} \label{eq:small_translations_small_injradius_base}
        \sum_{T \in \Gamma(0)} \mheat{t}{\tilde x}{T \tilde x}{\hyperbolic} \lesssim \sum_{\substack{k \neq 0 \\ \mdist{\tilde x}{\alpha^k \tilde x}{\hyperbolic} \leq 1}} \frac{1}{t} \exp \left( -\frac{\mdist{\tilde x}{\alpha^k \tilde x}{\hyperbolic}^2}{8t}\right) \, .
    \end{equation}
    Since $\mdist{\tilde x}{\alpha^k \tilde x}{\hyperbolic} \geq |k|\length{\gamma}$ always, the $k$ values in the sum in \eqref{eq:small_translations_small_injradius_base} satisfy $|k| \length{\gamma} \leq 1$. Denoting by $\rho$ the hyperbolic distance from $\tilde{x}$ to the imaginary axis, using \eqref{eq:injradius_estimate} and considering the Taylor expansion $\cosh(z) = 1+O(z^2)$ in \eqref{eq:distance_of_translates}, we find that 
    there exists a universal $C>0$ such that 
    \begin{equation*}
        \mdist{\tilde x}{\alpha^k \tilde x}{\hyperbolic} \geq \cosh^{-1}\left(1+C^2 \abs{k}^2 r^2 \right) \, .
    \end{equation*}
    Since $\cosh^{-1}(1+z) \asymp \sqrt{z}$ for $z\leq 2$ and $\cosh^{-1}(1+z) \asymp \log(z)$ for $z \geq 2$, there exists a universal constant $C>0$ such that
    \begin{equation*}
        \sum_{T \in \Gamma(0)} \mheat{t}{\tilde x}{T \tilde x}{\hyperbolic} 
        \lesssim \frac{1}{t} \sum_{k=1}^{\infty} \exp \left( -\frac{k^2 r^2 C^2}{8t}\right) + \frac{1}{t} \sum_{k=\ceil{2/r}+1}^{\infty} \exp \left( -\frac{C \log(k r)^2}{8t}\right) \, .
    \end{equation*}
    The first term on the right-hand side can be bounded as follows:
    \begin{align*}
        \frac{1}{t} \sum_{k=1}^{\infty} \exp \left( -\frac{k^2 r^2 C^2}{8t}\right) \nonumber        
        &\leq \frac{1}{t} \sum_{k=1}^{\infty} \exp \left( -\frac{(k^2+1) r^2 C^2}{16t}\right) \nonumber \\
        &\lesssim \frac{1}{t} \exp \left( -\frac{r^2 C^2}{16t}\right) \int_0^\infty  \exp \left( -\frac{y^2 r^2 C^2}{16t}\right)dy \nonumber\\    
        &\lesssim \frac{1}{t} \exp \left( -\frac{r^2 C^2}{16t}\right) \frac{\sqrt{t}}{r}  \, .        
    \end{align*}
    We may bound this last expression in two ways: first, bounding the exponential by $1$, we obtain 
    \begin{equation*}
        \frac{1}{t} \sum_{k=1}^{\infty} \exp \left( -\frac{k^2 r^2 C^2}{8t}\right) \lesssim \frac{1}{\sqrt{t}}\frac{1}{r} \, ,
    \end{equation*}
    which is bounded above by the right-hand side of \eqref{eq:ft_bound}; second, using the fact that $\frac{1}{x}\exp\left(-1/x\right) \leq 1$, we obtain 
    \begin{equation*}
        \frac{1}{t} \sum_{k=1}^{\infty} \exp \left( -\frac{k^2 r^2 C^2}{8t}\right) \lesssim \frac{\sqrt{t}}{r^3} \, ,
    \end{equation*}
    which is bounded above by the right-hand side of \eqref{eq:sqrt_time_bound}.
    The second term on the right-hand side can be bounded as follows:
    \begin{align}
        \frac{1}{t} \sum_{k=\ceil{2/r}+1}^{\infty} \exp \left( -\frac{C \log(k r)^2}{8t}\right) 
        &\lesssim \sum_{k=\ceil{2/r}+1}^{\infty} \frac{\log(kr)^2}{t}\exp \left( -\frac{C \log(k r)^2}{16t}\right) ^2  \nonumber \\
        &\lesssim \sum_{k=\ceil{2/r}+1}^{\infty} \exp \left( -\frac{C \log(k r)^2}{16t}\right) \nonumber \\
        &\lesssim \int_{2/r}^{\infty} \exp \left( -\frac{C \log(y r)^2}{16t}\right) ~dy \nonumber \\
        &\leq \frac{2}{r} \int_{0}^{\infty} \exp \left( -\frac{C z^2}{16t} + z \right) ~ dz \nonumber \\
        &\lesssim \frac{\sqrt{t}}{r} e^{16t} \label{eq:small_translations_small_injradius2} \, ,
    \end{align}    
    where the first inequality is because $kr \geq 2$, and the second is due to the fact that $\frac{1}{x}\exp(-1/x) \leq 1$. This is bounded by the right-hand side of both \eqref{eq:ft_bound} and \eqref{eq:sqrt_time_bound}, giving the desired result.
\end{proof}

\section{Proofs}
\begin{proof}[Proof of \cref{lem:ds_bounded_large_times}]
    Recall that by definition,
    \begin{equation*}
    \tracediff{t}{\S_k} - 1 = \int_{\S_k} \Big (\mheat{t}{x}{x}{\S_k} - \mheat{t}{0}{0}{\hyperbolic} - \frac{1}{\vol(\S_k)}\Big ) ~dx \, ,
    \end{equation*}
    and so
    \begin{align} \label{eq:splitting_S_large_t}
        \frac{\left| \tracediff{t}{\S_k} - 1 \right|}{t ~ \vol(\S_k)} &\leq \frac{1}{t ~ \vol(\S_k)}\int_{\S_k} \left( \mheat{t}{x}{x}{\S_k} - \frac{1}{\vol(\S_k)} \right) ~dx + \frac{\mheat{t}{0}{0}{\hyperbolic}}{t} \, .
    \end{align}
    By a well-known bound \cite[][Theorem 3.1]{davies_mandouvalos_heat_kernel_bounds} for the heat-kernel on the hyperbolic plane, we have $\mheat{t}{0}{0}{\hyperbolic} \lesssim e^{-t/4}/t$. The lemma then follows immediately from \cref{thm:convergence_to_uniform_large_t}.    
\end{proof}
~
\begin{proof}[Proof of \cref{lem:ds_uniformly_integrable_small_times}]
    Let $\gamma_1, \ldots, \gamma_s$ be the simple closed geodesics in $\S$ of length $\leq 2\sinh^{-1}(1)$, and denote their lengths by $\ell_1, \ldots, \ell_s$. We first give the following small-time bound
    \begin{equation}\label{eq:heat_kernel_trace_depends_on_local_geometry}    
        \int_{\S} \big[\mheat{t}{x}{x}{\S} - \mheat{t}{0}{0}{\hyperbolic} \big] ~ dx \lesssim \vol(\S)\sqrt{t} + \sum_{i=1}^{s} \min \Bigl\{ \frac{\sqrt{t}}{\ell_i^2}, \frac{1}{\sqrt{t}} \Bigr\} + \sum_{i \mid \ell_i \leq \sqrt{t}} \frac{1}{t^{1/4} ~ \sqrt{\ell_i}} \quad \forall t \in (0,1) \, .
    \end{equation}
    To show this, we partition the surface into a ``thin'' part and a ``thick'' part using \cref{lem:collar_lemma}. Let $C(\gamma_i)$ be the collars of the short geodesics. For a collar $C$, let $\truncated{C} := \{x \in C \mid \mathrm{dist}(x, \partial C) \geq 1\}$ be the ``truncated'' collar, obtained by cutting off the parts of the collar which are close to the boundary. Let $\thin{\S} = \union_{i=1}^{s}\truncated{C}(\gamma_i)$ be the union of all truncated collars and set $\thick{\S} := \S \backslash \thin{\S}$. The integral over the difference in heat-kernels can then be decomposed into
    		\begin{align*}
    			\displaystyle \int_{\S} \big [ \mheat{t}{x}{x}{\S} - \mheat{t}{0}{0}{\hyperbolic}\big ] ~dx = \int_{\thick{\S}} \big [\mheat{t}{x}{x}{\S} - \mheat{t}{0}{0}{\hyperbolic}\big ] ~dx + \int_{\thin{\S}} \big [\mheat{t}{x}{x}{\S} - \mheat{t}{0}{0}{\hyperbolic}\big ] ~dx \,.
    		\end{align*}
    The thick part is aptly named: for every $x \in \thick{\S}$, if $x \notin \union_i C_i$, then $\inj{x,\S} > \sinh^{-1}(1)$; while if $x \in C(\gamma_i) \backslash \truncated{C}(\gamma_i)$ and $\inj{x, \S} \leq \sinh^{-1}(1)$, then by \cref{lem:collar_lemma}, for $d = \mathrm{dist}(x, \partial C)$, 
    \begin{align*}
        \inj{x,\S_k} &= \sinh^{-1}\left(\cosh\left(\frac{1}{2}\length{\gamma_i}\right)\cosh(d) - \sinh(d) \right) \\
        &\geq \sinh^{-1}\left(\cosh(d) - \sinh(d) \right) \geq \sinh^{-1}\left(\cosh(1) - \sinh(1) \right) \, .
    \end{align*}
    Thus, the injectivity radius for points in the thick part is uniformly bounded below by some universal constant. By \eqref{eq:sqrt_time_bound}, we thus have
    \begin{equation*}  
        \mheat{t}{x}{x}{\S} - \mheat{t}{0}{0}{\hyperbolic} \lesssim \sqrt{t} \,
    \end{equation*}
    for all $0 < t \leq 1$ and every $x \in \thick{\S}$. This gives the first term in the right-hand side in the inequality \eqref{eq:heat_kernel_trace_depends_on_local_geometry}.
    
    For the thin part, no uniform pointwise estimate can be made on the heat-kernel. We instead integrate the difference $\mheat{t}{x}{x}{\S}- \mheat{t}{0}{0}{\hyperbolic}$ over every truncated collar $\truncated{C}(\gamma_i)$ separately. For brevity, denote $\ell := \length{\gamma_i}$, and for $x \in \S$, denote 
    \begin{equation*}
    r(x) := \inj{x,\S} \, .
    \end{equation*}
    We partition the truncated collar into two parts: points $x$ for which $r(x)^2 \geq t$, and points for which $r(x)^2 \leq t$.

    \paragraph{The case $\mathbf{r(x)^2 \geq t}$}
    For such points, we have
    \begin{align}
        \mheat{t}{x}{x}{\S}-\mheat{t}{0}{0}{\hyperbolic} ~&\substack{\text{\eqref{eq:sqrt_time_bound}} \\ \lesssim}~ \left(1+\frac{1}{r(x)^3} \right) \sqrt{t}  \nonumber \\
        &\substack{ \text{\eqref{eq:injradius_estimate}}  \\ \lesssim}~ \frac{\sqrt{t}}{\ell^3 \cosh(\rho)^3} \, , \label{eq:small_t_large_injradius_bound}
    \end{align}
    where $\rho$ is the distance of $x$ from the corresponding waist $\gamma_i$.
    Using \eqref{eq:injradius_estimate}, let $\beta_2 > 1$ be a universal constant such that $r(z)\leq \beta_2 \ell \cosh(\rho)$ for every $z \in \tilde{C}$. When $t < \beta_2^2 \ell^2$, we bound the integral of the heat kernel difference from above by integrating over all points in the truncated collar. Thus, for this range of $t$, 
    \begin{align}
        \int_{r(x)^2 \geq t} \big[ \mheat{t}{x}{x}{\S}-\mheat{t}{0}{0}{\hyperbolic} \big] ~dx &\leq \int_{0}^{2\pi} \int_{0}^{W(\gamma_i)-1} \big[ \mheat{t}{x}{x}{\S}-\mheat{t}{0}{0}{\hyperbolic} \big] \ell \cosh(\rho) ~d\rho d\theta \nonumber\\
        &\substack{\text{\eqref{eq:small_t_large_injradius_bound}} \\ \lesssim}~ \int_{0}^{\infty} \frac{\sqrt{t}}{\ell^2 \cosh(\rho)^2} d\rho = \frac{\sqrt{t}}{\ell^2} \tanh(\rho) \mid_0^\infty = \frac{\sqrt{t}}{\ell^2} \, .\label{eq:small_t_large_injradius_integral_bound1}
    \end{align}
    Now suppose that $t = \alpha^2 \beta_2^2 \ell^2$, where $\alpha \geq 1$. Since $r(x)^2 \geq t$, by the choice of $\beta_2$ we have $\cosh(\rho) \geq \alpha$. Thus, the integral can be bounded by integrating over all $\rho \geq \cosh^{-1}(\alpha)$:
    \begin{align*}
        \int_{r(x)^2 \geq t}\big[ \mheat{t}{x}{x}{\S}-\mheat{t}{0}{0}{\hyperbolic} \big] ~dx 
        &\leq \int_{0}^{2\pi} \int_{\cosh^{-1}\left( \alpha \right)}^{W(\gamma_i)-1} \big[ \mheat{t}{x}{x}{\S}-\mheat{t}{0}{0}{\hyperbolic} \big] \ell \cosh(\rho) ~d\rho d\theta \\    
        &\substack{\text{\eqref{eq:small_t_large_injradius_bound}} \\ \lesssim}~ \int_{\cosh^{-1}(\alpha)}^{\infty} \frac{\sqrt{t}}{\ell^2 \cosh(\rho)^2} ~ d\rho \\
        &= \frac{\sqrt{t}}{\ell^2} \left(1 - \tanh \left(\cosh^{-1}(\alpha) \right) \right) \leq \frac{\sqrt{t}}{\ell^2} \left(1 - \tanh \left(\log(\alpha) \right) \right) \, ,        
    \end{align*}    
    where the last inequality is due to the fact that $\cosh^{-1}(z) = \log \left(z + \sqrt{z^2-1}\right) \geq \log(z)$. Using the fact that $\tanh(z) = \frac{e^{2z}-1}{e^{2z}+1}$ and $\alpha = \sqrt{t}/\beta_2\ell$, we obtain
    \begin{equation} \label{eq:small_t_large_injradius_integral_bound2}
            \int_{r(x)^2 \geq t} \big[ \mheat{t}{x}{x}{\S}-\mheat{t}{0}{0}{\hyperbolic} \big] ~dx \lesssim \frac{\sqrt{t}}{\ell^2} \frac{1}{\alpha^2} \lesssim \frac{1}{\sqrt{t}}     \,
    \end{equation}
    for this range of $t$. Since $\sqrt{t}/\ell^2 \leq 1/\sqrt{t}$ for $t \leq \ell^2$, combining the bounds \eqref{eq:small_t_large_injradius_integral_bound1} and \eqref{eq:small_t_large_injradius_integral_bound2} gives the middle term in the right-hand side of equation \eqref{eq:heat_kernel_trace_depends_on_local_geometry}.

    \paragraph{The case $\mathbf{r(x)^2 \leq t}$}
    Since we are interested only in times $t \leq 1$, we use \cref{thm:li_yau} and obtain that
    \begin{equation}
        \mheat{t}{x}{x}{\S}-\mheat{t}{0}{0}{\hyperbolic} \leq \mheat{t}{x}{x}{\S}
        \lesssim \frac{1}{\vol(B(x,\sqrt{t}))} \, . \label{eq:beginning_of_li_yau}
    \end{equation}
    Bounding the integral $\int_{r(x)^2 \leq t}\big[ \mheat{t}{x}{x}{\S}-\mheat{t}{0}{0}{\hyperbolic} \big] ~dx$ is therefore a matter of giving a lower bound on the volume of the ball of radius $\sqrt{t}$ at $x$.
    \begin{figure}
        \centering
        \includegraphics[width=0.5 \textwidth]{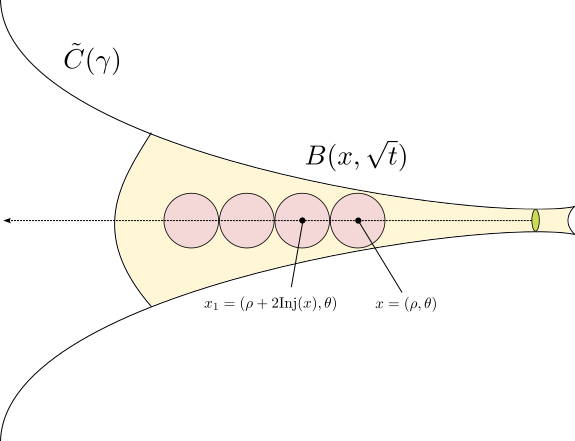}
        \caption{For large $t$, many copies of $B(x, r(x))$ can fit in $B(x,\sqrt{t})$.} 
        \label{fig:circles_on_collar}
    \end{figure}
    
    Let $m$ be the greatest integer such that $t \geq (2m+1)^2 r(x)^2$; since $t \geq r(x)^2$, we always have $m \geq 0$. Then the ball $B(x,\sqrt{t})$ contains $m+1$ disjoint balls of radius $r(x)$: If $x$ is given by the parameterization $(\rho,\theta)$, then there is one such ball at each point $x_j = (\rho + 2j r(x), \theta)$ for $j=0,\ldots, m$ (see Figure \ref{fig:circles_on_collar}). Note that since $t\leq 1$, by the truncation of $\truncated{C}(\gamma_i)$ from $C(\gamma_i)$, all the points of the form $(\rho+2j r(x),\theta)$ are still inside $C(\gamma_i)$ (even if they are not inside $\truncated{C}(\gamma_i)$), and so the injectivity radius is increasing in $j$; this ensures that injectivity radius at the center of each ball is at least $r(x)$. Since $\vol(B(x,s)) \geq s^2$, this gives
    \begin{equation*}
        \vol(B(x,\sqrt{t})) \geq (m+1) \cdot r(x)^2 \, .
    \end{equation*}
    This implies that when $t$ is written as $t = (2\alpha +1)^2 r(x)^2$ for real $\alpha > 0$, we have
    \begin{equation*}
        \vol(B(x,\sqrt{t})) \geq \frac {1}{2} (\alpha+1) \cdot r(x)^2 \geq \frac {1}{4} \sqrt{t} \cdot r(x) \, . 
    \end{equation*}
    If $x$ has cylindrical coordinates $(\rho,\theta)$ and satisfies $r(x)^2 \leq t$, then $\rho \leq \cosh^{-1}\left(\frac{\sqrt{t}}{\beta_1 \ell} \right)$, where $\beta_1$ is a universal constant, guaranteed to exist by \eqref{eq:injradius_estimate}, such that $r(z) \geq \beta_1 \ell\cosh(\rho)$ for all $z \in \tilde{C}$. Plugging the above volume estimate into the Li-Yau bound \eqref{eq:beginning_of_li_yau}, we thus obtain 
    \begin{align*}
        \int_{r(x)^2 \leq t} \big[ \mheat{t}{x}{x}{\S}-\mheat{t}{0}{0}{\hyperbolic} \big] ~dx &\lesssim \int_{r(x)^2 \leq t} \frac{1}{t^{1/2} r(x)} ~dx \\
        &\leq \int_0^{2\pi} \int_{0}^{\cosh^{-1}\left(\frac{\sqrt{t}}{\beta_1 \ell} \right)} \frac{1}{t^{1/2} r(x)}\ell \cosh(\rho) ~d\rho d\theta \\
        &\lesssim \int_{0}^{\cosh^{-1}\left(\frac{\sqrt{t}}{\beta_1 \ell} \right)} \frac{1}{t^{1/2}} ~d\rho  \\
        &= \cosh^{-1}\left(\frac{\sqrt{t}}{\beta_1 \ell} \right) \frac{1}{t^{1/2}} \lesssim \frac{1}{\sqrt{\ell} \cdot t^{1/4}}\, ,    
    \end{align*}
    where the last inequality is due to the fact that $\cosh^{-1}(z) \leq \sqrt{2z}$. Noticing that the inequality $r(x)^2 \leq t$ can only hold when $t \geq \ell^2$, this gives the third and last term in the right-hand side of
    equation \eqref{eq:heat_kernel_trace_depends_on_local_geometry}.

    With \eqref{eq:heat_kernel_trace_depends_on_local_geometry} at hand, we now show that if $\S_k$ have uniformly integrable geodesics, then $\frac{\tracediff{t}{\S_k}}{t\vol(\S_k)}$ is uniformly integrable on $(0,1)$. By \eqref{eq:heat_kernel_trace_depends_on_local_geometry},
    \begin{equation*}
    \frac{\tracediff{t}{\S_k}}{t ~\vol(\S_k)} \lesssim 
        \left(\frac{1}{\sqrt{t}} 
        + \frac{1}{\vol(\S_k)}\sum_{i=1}^{s} \min \Bigl\{ \frac{1}{\sqrt{t}~\ell_i^2}, \frac{1}{t^{3/2}}\Bigr\} 
        + \frac{1}{\vol(\S_k)}\sum_{i \mid \ell_i \leq \sqrt{t}} \frac{1}{t^{5/4} ~ \sqrt{\ell_i}}\right) := \gsmall(t) \, .
    \end{equation*}
    To show that $\gsmall$ is indeed a uniformly integrable family, note first that $\gsmall(t)$ are all decreasing in $t$, and so it suffices to show that their $\mathrm{L}^1$ norms are all bounded, and that for every $\eps > 0$ there exists a $\delta > 0$ so that for all $k$ large enough,
    \begin{equation*}
        \int_0^\delta \gsmall(t) ~dt < \eps \, .
    \end{equation*}
    Let us now bound the integrals of the three terms that make up $\gsmall$.
    
    The first term, $\frac{1}{\sqrt{t}}$, is integrable.
    
    For the second term, for every $\delta > 0$ and every $i=1,\ldots, s$, we have
    \begin{align}
        \int_0^\delta \min \Bigl\{ \frac{1}{\sqrt{t}~\ell_i^2}, \frac{1}{t^{3/2}}\Bigr\}dt &= \int_0^{\min(\ell_i^2, \delta)}\frac{1}{\sqrt{t}~\ell_i^2}~dt + \int_{\min(\ell_i^2,\delta)}^\delta \frac{1}{t^{3/2}}~dt = 2\frac{\sqrt{t}}{\ell_i^2} \,\,\Big \vert_0^{\min \{\ell_i^2,\delta\}} - \frac{2}{\sqrt{t}} \,\, \Big \vert_{\min\{\ell_i^2, \delta\}}^\delta  
        \nonumber \\
        &= 2 \left(\min \Bigl\{ \frac{1}{\ell_i}, \frac{\sqrt{\delta}}{\ell_i^2}\Bigr\} - \frac{1}{\sqrt{\delta}} + \frac{1}{\min \{\ell_i, \sqrt{\delta}\}} \right) \, .    \label{eq:gsmall_integral}
    \end{align}
    Let $\eta > 0$. Summing over all $i=1,\ldots, s$ and separating the sum into $\ell_i \leq \eta$ and $\ell_i \geq \eta$, we get
    \begin{align*}
        \frac{1}{\vol(\S_k)} \sum_1^s\int_0^\delta \min \Bigl\{ \frac{1}{\sqrt{t}~\ell_i^2}, \frac{1}{t^{3/2}}\Bigr\}dt &\leq \frac{2}{\vol(\S_k)}\left( \sum_{\ell_i < \eta} \frac{1}{\ell_i} + \sum_{\ell_i \geq \eta} \frac{\sqrt{\delta}}{\ell_i^2 } + \sum_{\ell_i < \sqrt{\delta}} \frac{1}{\ell_i} \right) \\    
        &\leq 2\invlengths_\eta(\S_k) + \frac{2}{\vol(\S_k)}\sum_{\ell_i \geq \eta} \frac{\sqrt{\delta}}{\eta^2} + 2\invlengths_{\sqrt{\delta}}(\S_k) \, .
    \end{align*}
    First choosing $\eta$ small enough so that $\invlengths_\eta < \eps$ for large enough $k$ and then choosing $\delta = \eps \eta^2$ gives that $\int_0^\delta \gsmall (t)dt \lesssim \eps$ for all large $k$. Note that for $\delta  = 1$, \eqref{eq:gsmall_integral} yields $\int_0^1 \gsmall(t)dt \leq 4\invlengths_{2\sinh^{-1}(1)}(\S_k)$, which, by our assumption on $\invlengths_\eta(\S_k)$, is bounded independently of $k$. 
    
    For the third term, we have
    \begin{align*}
        \int_0^\delta \sum_{i \mid \ell_i \leq \sqrt{t}} \frac{1}{t^{5/4} ~ \sqrt{\ell_i}} dt &= \sum_{i=1}^{s} \int_{\min(\ell_i^2, \delta)}^\delta \frac{1}{t^{5/4}~\sqrt{\ell_i}} dt = -4\sum_{i=1}^{s} \frac{1}{t^{1/4} \sqrt{\ell_i}}\,\,\Big \vert_{\min(\ell_i^2, \delta)}^\delta 
        \\
        &= \sum_{i=1}^{s} \left(\frac{4}{\min(\ell_i,\delta^{1/4}\sqrt{\ell_i})} - \frac{4}{\delta^{1/4}\sqrt{\ell_i}}\right) \leq \sum_{\ell_i < \sqrt{\delta}} \frac{4}{\ell_i} = 4\invlengths_{\sqrt{\delta}} \, .
    \end{align*}
    Taking $\delta = 1$ gives that this term is uniformly bounded, while taking $\delta$ small enough so that $\invlengths_{\sqrt{\delta}} < \eps$ for $k$ large enough gives uniform integrability of this term. This shows that if the $\S_k$ have uniformly integrable geodesics, then $\frac{\tracediff{t}{\S_k}}{t \vol(\S_k)}$ are uniformly integrable.

    For the converse, it suffices to show that for every $\eta <1/2$ and every compact hyperbolic surface $\S$,
    \begin{equation} \label{eq:lower_bound_on_ds}
        \int_{0}^{\eta} \frac{\tracediff{t}{\S}}{t~\vol(\S)} ~dt \gtrsim \invlengths_{\eta}(\S) \, .        
    \end{equation}
    Indeed, since $\S_k$ do not have uniformly integrable geodesics, there exists an $\eps > 0$ such that for every $\eta > 0$, $\invlengths_{\eta}(\S_k) >\eps$ for infinitely many $k$. The above inequality then ensures that for every $\eta >0$, $\int_{0}^{\eta} \frac{\tracediff{t}{\S_k}}{t~\vol(\S_k)} ~dt \gtrsim \eps$ for infinitely many $k$. To prove this inequality, let $\gamma$ be a simple closed geodesic of length $\ell \leq \eta$, let $C$ be its collar, and let $\alpha$ be the group element which generates $C$. For every $x \in C$, we have 
    \begin{equation}
        \mheat{t}{x}{x}{\S} - \mheat{t}{0}{0}{\hyperbolic} =\sum_{T \neq \mathrm{Id} , T \in \Gamma} \mheat{t}{\tilde x}{T \tilde x}{\hyperbolic} \geq \mheat{t}{\tilde x}{\alpha \tilde x}{\hyperbolic} \, ,
    \end{equation}
    where $\S = \Gamma \backslash \hyperbolic$, $\tilde{x}$ is the inverse image of $x$ under to projection to $\S$. Let $A = \{x \in \S \mid \dist{x}{\gamma} \leq 1\}$ be the set of points at distance at most $1$ from $\gamma$. Since $\ell < 1/2$, we have $W(\gamma) > 2$, and so $A \subseteq C$ and $\dist{A}{\partial C} > 1$.

    As is shown in the proof of \cref{lem:counting_argument_by_injradius}, for any nontrivial $T \neq \alpha^k$, $\mdist{\tilde{x}}{T\tilde{x}}{\hyperbolic} \geq \dist{x}{\partial C} > 1$ for all $x \in A$. On the other hand, $\mdist{\tilde{x}}{\alpha\tilde{x}}{\hyperbolic} < 1$ for $\ell < 1/2$, as calculated by \eqref{eq:alternate_injectivity_in_collar} (in an infinite hyperbolic collar, the injectivity radius of $z$ is equal to $\frac{1}{2}\mdist{z}{\alpha z}{\hyperbolic}$). Thus $\inj{x} = \frac{1}{2}\mdist{\tilde{x}}{\alpha\tilde{x}}{\hyperbolic}$. By \eqref{eq:injradius_estimate}, $\inj{x} \asymp \ell$, and so 
    \begin{equation} \label{eq:group_element_cant_move_far}
        \mdist{\tilde{x}}{\alpha\tilde{x}}{\hyperbolic} \asymp \ell \, .
    \end{equation}
    By \cite[][Theorem 3.1]{davies_mandouvalos_heat_kernel_bounds}, the heat-kernel in the hyperbolic plane is proportional to 
\begin{equation*}
    \mheat{t}{x}{y}{\hyperbolic} \sim \frac{1}{t}\exp\left( -\frac{1}{4}t - \frac{\mdist{x}{y}{\hyperbolic}^2}{4t} - \frac{1}{2}\mdist{x}{y}{\hyperbolic} \right) \frac{1+\mdist{x}{y}{\hyperbolic}}{1 + \mdist{x}{y}{\hyperbolic} + t} \, .
\end{equation*}
    Together with \eqref{eq:group_element_cant_move_far}, this implies that for $x \in A$ and for $t \in [\frac{1}{2} \ell^2, \ell^2]$,
    \begin{equation*}
        \mheat{t}{\tilde x}{\alpha \tilde x}{\hyperbolic} \gtrsim \frac{1}{t} \, . 
    \end{equation*}
    We then have 
    \begin{align*}
        \int_0^\eta \frac{\int_{C} [\mheat{t}{x}{x}{\S} - \mheat{t}{0}{0}{\hyperbolic}]~dx}{t}dt 
        &\gtrsim \int_{\frac{1}{2}\ell^2}^{\ell^2} \int_{A} \frac{1}{t^2} ~dx dt \\
        &\asymp \int_{\frac{1}{2}\ell^2}^{\ell^2} \int_{0}^{1} \frac{1}{t^2} \ell \cosh(\rho) ~d\rho dt \asymp \frac{1}{\ell} \, .    
    \end{align*}
    Summing over all $\gamma$ of length $\leq \eta$ and dividing by $\vol(\S)$ gives the desired result.

\end{proof}

\begin{proof} [Proof of \cref{thm:logdet_main_theorem}]
Observe that if $\S_k$ has uniformly integrable geodesics, then so does every subsequence $\S_{k_m}$. Thus, in order to show item \eqref{enu:main_theorem_uniform_integrable} for the original sequence $\S_k$, it suffices to show it for every subsequence such that $\frac{\log \det \Delta_{\S_{k_m}}}{\vol(\S_{k_m})}$ converges. On the other hand, if $\S_k$ does not have uniformly integrable geodesics, then different subsequences may or may not have uniformly integrable geodesics. Still, item \eqref{enu:main_theorem_non_uniform_integrable} will follow for $\S_k$ if we can show that for every subsequence without uniformly integrable geodesics such that $\frac{\log \det \Delta_{\S_{k_m}}}{\vol(\S_{k_m})}$ converges, the limit is strictly smaller than $E_\mu$ if it converges to a finite number, and is $-\infty$ if not. Dropping the subindex, we can therefore restrict ourselves to subsequences $\S_k$ such that $\frac{\log \det \Delta_{\S_k}}{\vol(\S_k)}$ converges.

Assume first that $\sup_k \invlengths_1(\S_k) = \infty$ (in particular, the $\S_k$ do not have uniformly integrable geodesics). Then by \eqref{eq:the_starting_point}, by \cref{lem:ds_bounded_large_times} and by \eqref{eq:lower_bound_on_ds} from the proof of \cref{lem:ds_uniformly_integrable_small_times}, we have
\begin{align*}
    \inf_{k} \frac{\log \det \Delta_{\S_k}}{\vol(\S_k)} &= E_\hyperbolic  +\inf_{k}\left(\frac{\gamma_0}{\vol(\S_k)} - \int_0^1\frac{\tracediff{t}{\S_k}}{t ~ \vol(\S_k)}dt - \int_1^\infty \frac{\tracediff{t}{\S_k}-1}{t ~ \vol(\S_k)} dt \right) \\
    &\lesssim E_\hyperbolic + \gamma_0 - \sup_k \left( \invlengths_1(\S_k) - \sqrt{1 + \invlengths_1(\S_k)} \right) \, .
\end{align*}
Since $\sup_k \invlengths_1(\S_k) = \infty$, the right-hand side of the above inequality is equal to $-\infty$. 
Since the sequence $\frac{\log \det \Delta_{\S_k}}{\vol(\S_k)}$ converges and the infimum is $-\infty$ we get
\begin{align*}
    \lim_{k \to \infty} \frac{\log \det \Delta_{\S_k}}{\vol(\S_k)} = \liminf_{k \to \infty} \frac{\log \det \Delta_{\S_k}}{\vol(\S_k)} = \inf_{k} \frac{\log \det \Delta_{\S_k}}{\vol(\S_k)} = -\infty \, ,
\end{align*}
meaning the limit of the sequence is not finite, satisfying item \eqref{enu:main_theorem_non_uniform_integrable} of \cref{thm:logdet_main_theorem}.

We now assume that $\sup_k \invlengths_1(\S_k) := M < \infty$. This does not imply that $\S_k$ has uniformly integrable geodesics; however it does imply that the logarithms of the reciprocal lengths of the geodesics must be uniformly integrable. Indeed, let $\eps >0$. Since $x\log(1/x) \to 0$ as $x \to 0$, there exists an $\eta >0$ such that $x\log(1/x) < \eps/M$ for all $x < \eta$, implying that for every $k$,
\begin{equation} \label{eq:log_uniformity}
    \frac{1}{\vol(\S_k)}\sum_{\ell_i < \eta}\log\left(1/\ell_i\right) \leq \frac{1}{\vol(\S_k)}\sum_{\ell_i < \eta}\frac{\eps}{M}\frac{1}{\ell_i} \leq \eps \, .
\end{equation}
Under this condition, we now show that for every fixed $t > 0$, the limit $\lim_{k\to\infty} \frac{\tracediff{t}{\S_k}}{\vol(\S_k)}$ exists. 

Recall that $\mathcal{H}$ is the space of all rooted hyperbolic surfaces, endowed with the Gromov-Hausdorff metric. Fix $t > 0$, and let $h_t:\mathcal{H} \to \reals$ be given by $h_t(\S,x) = \mheat{t}{x}{x}{\S} - \mheat{t}{0}{0}{\hyperbolic}$. The function $h_t$ is continuous under the Gromov-Hausdorff topology (see e.g. \cite[Theorem 2.6, Theorem 5.54, and Theorem 5.59]{ding_heat_kernels_and_greens_functions}. These theorems show that $h_t$ converges on a sequence of  compact surfaces (respectively non-compact) which converge to a target compact (respectively non-compact) limiting surface. The same methods work also for the remaining two cases.). Let $\mu_k$ be the uniform measure on the set $\{(\S_k, x) \mid x \in \S_k \} \subseteq\mathcal{H}$, and define $X_k = h_t(\S_k,x)$, where $(\S_k,x) \sim \mu_k$ and $X = h_t(\mathbf{\S}_\infty, \mathbf{x}_0)$, where $(\mathbf{\S}_\infty, \mathbf{x}_0)$ distributes as $\mu$. Showing that the limit $\lim_{k\to\infty} \frac{\tracediff{t}{\S_k}}{\vol(\S_k)}$ exists is then equivalent to showing that $\lim_{k\to \infty}\mathbb{E}[X_k]$ exists. 
If $h_t$ were bounded, convergence of $\mathbb{E}[X_k] = \int h_t ~ d\mu_k$ would have followed from the weak convergence of $\mu_k$ to $\mu$. However, $h_t$ is unbounded, and we will utilize the uniform integrability of the logarithms, as in equation \eqref{eq:log_uniformity}. We will use Vitali's convergence theorem for random variables (\cite[][Theorem 25.12]{billinglsey_probability_and_measure}), which states that if $X_k$ are a sequence of uniformly integrable random variables which converge in distribution to $X$, then $X$ is integrable and $\mathbb{E}[X_k] \to \mathbb{E}[X]$. Here, a sequence of random variables is uniformly integrable if for every $\eps >0$ there is an $M > 0$ such that $\sup_k \mathbb{E}[X_k \one_{X_k \geq M}] < \eps$.

First note that the $X_k$ converge in distribution to $X$: From the continuous mapping theorem \cite[][Corollary 1 of Theorem 25.7]{billinglsey_probability_and_measure}, it follows that if a sequence of random variables $Y_k \in \mathcal{H}$ converge in distribution to a random variable $Y \in \mathcal{H}$, and $h: \mathcal{H} \to \reals$ is continuous, then $h(Y_n)$ converges in distribution to $h(Y)$. Letting $Y_k \sim \mu_k$, $Y \sim \mu$ and taking $h = h_t$, we get that $X_k$ converges to $X$ in distribution. 

It remains to show that the $X_k$ are uniformly integrable. By \eqref{eq:ft_bound}, there exists a function $f(t)$ such that $h_t(\S,x) < f(t)\left(1 + 1/\inj{x}\right)$. It therefore suffices to show that for every $\eps > 0$, there is a $\delta > 0$ so that $\int \frac{1}{\inj{x}} \one_{\inj{x} < \delta}d\mu_k \leq \eps$ for every $k$. A calculation shows that in a collar $C$ with waist of length $\ell$,  
\begin{align*}
    \int_C \frac{1}{\inj{x}}~dx \lesssim \int_0^{W(\ell)} \int_0^{2\pi} \frac{1}{\ell \cosh(\rho)} \ell \cosh(\rho)~d\theta d\rho \asymp \log \left(1/\ell \right)  \, .
\end{align*}
Thus, by \eqref{eq:log_uniformity}, there exists a $\delta < \sinh^{-1}(1)$ small enough so that if $\ell_1, \ldots, \ell_s$ are the lengths of simple closed geodesics shorter than $\delta$ in $\S_k$, then
\begin{equation} \label{eq:log_uniformity_suffices}
    \int_{\S_k} \frac{\one_{\inj{x} < \delta}}{\inj{x}} ~d\mu_k \lesssim \frac{1}{\vol(\S_k)} \sum_{i=1}^{s} \log(1/\ell_i) \leq \eps \, .
\end{equation}
Applying Vitali's theorem shows that $\lim_{k\to\infty} \frac{\tracediff{t}{\S_k}}{\vol(\S_k)}$ exists and is equal to $\int \big[\mheat{t}{x_0}{x_0}{\S_\infty} - \mheat{t}{0}{0}{\hyperbolic} \big] ~ d\mu$.

Having shown this, consider the limit of the large-times integral, $\lim_{k\to\infty} \int_1^\infty\frac{\tracediff{t}{\S_k}-1}{t ~ \vol(\S_k)}dt$. By \cref{lem:ds_bounded_large_times} and the fact that $\sup_k \invlengths_1(\S_k) < \infty$, we have that $\frac{\abs{\tracediff{t}{\S_k}-1}}{t~\vol(\S_k)}$ is uniformly bounded by an integrable function, and so by the dominated convergence theorem,
\begin{equation} \label{eq:exchange_large_times}
    \lim_{k\to\infty} \int_1^\infty\frac{\tracediff{t}{\S_k}-1}{t ~ \vol(\S_k)}dt = \int_1^\infty \lim_{k\to\infty}  \frac{\tracediff{t}{\S_k}-1}{t ~ \vol(\S_k)}dt = \int_1^\infty  \frac{\int \big[\mheat{t}{x_0}{x_0}{\S_\infty} - \mheat{t}{0}{0}{\hyperbolic} \big] ~ d\mu}{t} ~ dt \, .
\end{equation}

For the small-times integral, we use Vitali's convergence theorem for functions \cite[][Theorem 16.14]{billinglsey_probability_and_measure}, which states that if ${f_k}$ is a family of non-negative functions which converges almost everywhere to a measurable function $f$, then $f_k$ are uniformly integrable if and only if $f_k \xrightarrow{L^1} f$. Together with \cref{lem:ds_uniformly_integrable_small_times}, this implies that 
\begin{equation} \label{eq:exchange_small_times}
    \lim_{k\to\infty} \int_0^1\frac{\tracediff{t}{\S_k}}{t ~ \vol(\S_k)}dt = \int_0^1 \frac{\int \big[\mheat{t}{x_0}{x_0}{\S_\infty} - \mheat{t}{0}{0}{\hyperbolic} \big] ~ d\mu}{t} ~ dt \, 
\end{equation}
if and only if $\S_k$ has uniformly integrable geodesics. When $\S_k$ has uniformly integrable geodesics, \eqref{eq:the_starting_point}, \eqref{eq:exchange_large_times} and \eqref{eq:exchange_small_times} give item \eqref{enu:main_theorem_uniform_integrable} in \cref{thm:logdet_main_theorem}. When $\S_k$ does not have uniformly integrable geodesics, Fatou's lemma and \eqref{eq:exchange_large_times} imply that 
\begin{align*}
        \limsup_{k \to \infty} \frac{\log \det \Delta_{\S_k}}{\vol(\S_k)} \leq E_\mu \, ;
\end{align*}
since the limit exists, by \eqref{eq:exchange_small_times} equality cannot occur, and we obtain item \eqref{enu:main_theorem_non_uniform_integrable} in \cref{thm:logdet_main_theorem}.
\end{proof} 
~
\begin{proof}[Proof of \cref{cor:convergence_to_plane}]
    Since $\hyperbolic$ is the universal cover of any hyperbolic surface, we have $\mheat{t}{x}{x}{\S} > \mheat{t}{0}{0}{\hyperbolic}$ for all $\S \neq \hyperbolic$ and all $t \in (0,\infty)$. Thus, if $\lim_{k\to\infty} \frac{\log \det \Delta_{\S_k}}{\vol(\S_k)} = E_{\hyperbolic}$, then $\mheat{t}{x}{x}{\S} = \mheat{t}{0}{0}{\hyperbolic}$ for all $t$ $\mu$-a.s, implying that $\mathbf{\S_\infty} = \hyperbolic$ $\mu$-a.s.
\end{proof}

\section*{Acknowledgements} The authors are supported by ERC consolidator grant 101001124 (UniversalMap) as well as ISF grants 1294/19 and 898/23.

\printbibliography
\end{document}